\documentclass[11pt]{article}
\usepackage{amsmath,amsthm,bbm,}

\newcommand{\R}{\mathbbm{R}}

\newcommand{\nothing}[1]{}
\newtheorem{theorem}{Theorem}[section]

\newtheorem{open problem}{Open problem}[section]

\newcommand{\ou}{{\overline{u}}}
\newcommand{\ot}{{\overline{t}}}

\newcommand{\De}{{\Delta}}
\newcommand{\na}{{\nabla}}
\newcommand{\pa}{{\partial}}
\newcommand{\cn}{{\nu}}

\newcommand{\om}{{\omega}}
\newcommand{\Om}{{\Omega}}
\newcommand{\Ga}{{\Gamma}}

\newcommand{\m}{{\mu}}

\begin{document}
\title{\bf On the extension to slip boundary conditions of a Bae and Choe regularity criterion for
the Navier-Stokes equations. The half-space case.}
\author{H.~Beir{\~a}o~da Veiga \footnote{Partially supported  by FCT (Portugal) under grant UID/MAT/04561/3013.},\\
Department of Mathematics,\\
Pisa University, Italy.\\
email: bveiga@dma.unipi.it} \maketitle

\vspace{0.2cm}

\begin{abstract}
This notes concern the sufficient condition for regularity of
solutions to the evolution Navier-Stokes equations known in the
literature as Prodi-Serrin's condition. H.-O. Bae and H.-J. Choe
proved in a 1999 paper that, in the whole space $\,\R^3\,$, it is
merely sufficient that two components of the velocity satisfy the
above condition. Below, we extend the result to the half-space case
$\,\R^n_+\,$ under slip boundary conditions. We show that it is
sufficient that the velocity component \emph{parallel} to the
boundary enjoys the above condition. Flat boundary geometry seems
not essential, as suggested by some preliminary calculations in
cylindrical domains.
\end{abstract}
\section{Introduction.}\label{uns}
These notes concern sufficient conditions for regularity of
solutions to the evolution Navier-Stokes equations related to the so
called Prodi-Serrin's condition, see \eqref{ps}. In reference
\cite{baechoe} the authors proved, in the whole space case, that it
is sufficient that two components of the velocity satisfy the above
condition. Below we extend this result to the half-space case
$\,\R^n_+\,$ under slip boundary conditions. The structure of the
proof follows Bae and Choe paper, by adding a suitable control of
some boundary integrals which, clearly, were not present in the
whole space case.\par%
We assume that readers are acquainted with the main literature on
the subject. In particular, we will not repeat well know
notation as, for instance, Sobolev spaces notation, and so on.%

\vspace{0.2cm}

In the sequel we are interested in the evolution Navier-Stokes
equations in the half-space $\,\R^n_+ =\,\{x:\,x_n>\,0\,\}\,,$
$n\geq\,3\,,$
\begin{equation}\label{nse-2}
\left\{\begin{aligned}%
&\pa_t\,u+(u\cdot\nabla)\,u-\,\m\,\triangle
u+\nabla \pi =0\,,\\
&\nabla\cdot u=0 \quad\text{in }\R^n_+ \times (0,T]\,; \\
&u(x,0)=\,u_0(x)\qquad\text{in }\R^n_+ \,,
\end{aligned}\right.
\end{equation}
under the classical Navier slip boundary conditions without
friction. See \cite{navier} and \cite{Ser59}. On flat portions of
the boundary this condition reads
\begin{equation}%
\label{bboundary} \left\{
\begin{split}
&u_n=0,\\
&\pa_n\, u_j=0,\quad 1\leq j\leq n-\,1 \,.\
\end{split}
\right.
\end{equation}
In the half-space case we will use this formulation. Let us recall
that, for $n=\,3\,,$ the above slip boundary condition may also be
written in the form $\,u_n=0\,,$ plus $\,\om_j=\,0,$ for
$\,j=\,1,\,2\,,$
where $\omega=\nabla\times u\,$ is the vorticity field.\par%
It is well know that weak solutions $\,u\,$ satisfying the so called
Prodi-Serrin's condition
\begin{equation}
\label{ps}%
u \in\, L^q(0,\,T;\,L^p(\R^n_+)\,)\,, \quad \frac2q
+\,\frac{n}p\leq\,1\,, \quad\, p>\,n
\end{equation}
are strong, namely
\begin{equation}
\label{regil}%
u \in\,L^{\infty}(0,\,T;\,H^1(\R^n_+)
\cap\,L^2(0,\,T;\,H^2(\R^n_+)\,)\,.
\end{equation}
The proof is classical. Furthermore, strong
solutions are smooth, if data and domain are also smooth.\par%
It is well known that the above results hold in a very large class
of domains $\,\Omega\,,$ under suitable boundary conditions. We
assume this kind of results well known to the reader. In particular,
the result is well known in the whole space $\,\R^n\,,$ which is our
departure point. In fact, consider the Navier-Stokes equations
\eqref{nse-2} with $\R^n_+\,$ replaced by $\,\R^n\,.$
Differentiating both sides of the first equation \eqref{nse-2} with
respect to $x_k$, taking the scalar product with $\pa_k \,u\,$,
adding over $k\,,$ and integrating by parts over $\R^n$, one shows
that
\begin{equation}
\label{estimun}%
\frac12\,\frac{d}{dt}\,\int_{\R^n}\,|\na\,u|^2 \,dx +\,\m\,
\int_{\R^n}\,|\na^2\,u |^2\,dx =\,-\,\int_{\R^n} \,\na
\,[(u\cdot\nabla)\,u\,]\cdot\na\,u  \, dx\,,
\end{equation}
where obvious integrations by parts have been done. Clearly, no
boundary integrals appear. A last integration by parts shows that
\begin{equation}
\label{canja}%
\big|\,\int_{\R^n} \,\na \,[(u\cdot\nabla)\,u\,]\cdot\na\,u\, dx
\big|\leq\, c(n)\, \int_{\R^n} \, |\,u\,|\,|\na\,u |\,|\na^2\,u |\,
dx \,.
\end{equation}
So
\begin{equation}
\label{estimdos}%
\frac12\,\frac{d}{dt}\,\int_{\R^n}\,|\na\,u|^2 \,dx +\,\m\,
\int_{\R^n}\,|\na^2\,u |^2\,dx \leq\, c(n)\, \int_{\R^n} \,
|\,u\,|\,|\na\,u |\,|\na^2\,u |\, dx
\end{equation}
follows.\par%
By appealing to the Prodi-Serrin's assumption \eqref{ps} applied to
the term $ |\,u\,|\,$ present in the right hand side of
\eqref{estimdos}, well know devices lead to the desired regularity
result \eqref{regil} in $\,\R^n \,$ (these same devices may be seen
in section \ref{dos}, in connection with the similar estimate
\eqref{estibom}). In the proof, the crucial property is that the
term $ |\,u\,|\,$ in the right hand side of estimate \eqref{canja}
(hence also in \eqref{estimdos}) enjoys the Prodi-Serrin's
condition. In reference \cite{baechoe}, see also \cite{baechoe-2},
H.-O.~Bae and H.-J- Choe succeed in replacing, in the right hand
side of \eqref{canja}, the term $\,|u|\,$ simply by $\,|\ou|\,$,
where $\,\ou\,$ is an arbitrary $n-1$ dimensional component
\begin{equation}
\label{serves}%
\ou=\,(u_1,\,...,\,u_{n-1},\,0)
\end{equation}
of the velocity $\,u\,.$ In other words, they succeed in improving
the quite obvious estimate \eqref{canja}, by showing the much
stronger estimate
\begin{equation}
\label{bacho}%
\big|\,\int_{\R^n} \,\na \,[(u\cdot\nabla)\,u\,]\cdot\na\,u\, dx
|\leq\, c(n)\, \int_{\R^n} \, |\ou|\,|\na\,u |\,|\na^2\,u |\, dx\,.
\end{equation}
Hence the estimate \eqref{estimdos} holds with $\,|\,u|\,$ replaced
by $\,|\,\ou|\,.$ The classical $\,|\,u|-$proof applies as well
after this substitution. In this way the authors proved that
\eqref{regil} holds if merely $\,\ou\,$ (instead of $\,u\,$)
satisfies the Prodi-Serrin's condition. A quite unexpected result,
at that time, may be not yet sufficiently exploited. Clearly, in the
whole space case, $\,\ou\,$ may be any $\,n-1\,$ dimensional
component of the velocity.\par%
The proof of the estimate \eqref{bacho} is based on a clever
analysis of the structure of the integral on the left hand side of
this equation.
\vspace{0.2cm}

The first aim of these notes is to prove equation \eqref{bacho} in
the half space $\,\R^n_+\,$
\begin{equation}
\label{bacho-2}%
\big|\,\int_{\R^n_+} \,\na \,[(u\cdot\nabla)\,u\,]\cdot\na\,u\, dx
|\leq\, c(n)\, \int_{\R^n_+} \, |\ou|\,|\na\,u |\,|\na^2\,u |\,
dx\,,
\end{equation}
under slip boundary conditions. As a consequence, the estimate
\eqref{estimdos} holds with $\,|u|\,$ replaced by $\,|\ou|\,$ and
$\,\R^n\,$ replaced by $\,\R^n_+\,$. It readily follows, as in the
classical case, that solutions to the above boundary value problem
are regular provided that $\,\ou\,$ satisfies the Prodi-Serrin's
condition \eqref{ps}.
\begin{theorem}
\label{theocas}%
Let $\,u\,$ be a solution to the Navier-Stokes equations
\eqref{nse-2} in $\,\R^n_+\,$ under the the slip boundary conditions
\eqref{bboundary}. Furthermore, let  $\,\ou\,$ be the parallel to
the boundary component of the velocity $\,u\,,$ given by
\eqref{serves}. If
\begin{equation}
\label{bachol}%
\ou \in\, L^q(0,\,T;\,L^p(\R^n_+)\,)\,, \quad \frac2q
+\,\frac{n}p\leq\,1\,, \quad\, p>\,n\,,
\end{equation}
then \eqref{regil} holds.
\end{theorem}
Alternatively, the proof of the above result could be done by
appealing to a reflection principle, see \cite{bvcrgr}. However this
does not help possible extensions to non flat boundaries. In this
direction, in a forthcoming paper, we will consider the cylindrical
coordinates case.

\vspace{0.2cm}

We end this section by quoting the very recent paper \cite{baewolf}
where the authors proved the local, interior, regularizing effect of
the Prodi-Serrin's condition only on two velocity components. It
would be of interest to extend this result to arbitrary, smooth,
coordinates (orthogonal for instance).
\section{Extension to boundary value problems}\label{dos}
In this section we prove equation \eqref{bacho-2}. Our approach adds
to that followed in reference \cite{baechoe} an accurate control of
the boundary integrals, clearly not present in the whole space case.
To obtain the explicit form of these integrals, we have to turn
back to the volume integrals.\par%
For notational convenience we set
$$
\Ga=\,\{ \,x:\, x_n =\,0\,\}\,.
$$
The first steep is to prove \eqref{estimun}, now with $\,\R^n\,$
replaced by $\,\R^n_+\,,$ namely
\begin{equation}
\label{estimum-2}%
\frac12\,\frac{d}{dt}\,\int_{\R^n_+}\,|\na\,u|^2 \,dx +\,\m\,
\int_{\R^n_+}\,|\na^2\,u |^2\,dx =\,-\,\int_{\R^n_+} \,\na
\,[(u\cdot\nabla)\,u\,]\cdot\na\,u  \, dx\,.
\end{equation}
By following the $\,\R^n\,$ case, we differentiate both sides of the
first equation \eqref{nse-2} with respect to $x_k$, take the scalar
product with $\pa_k \,u\,$, add over $k\,,$ and integrate by parts
over $\R^n_+$. Now additional boundary integrals appear. We start
from the $\,\triangle u\,$ term. One has
$$
-\,\int_{\R^n_+}\, \na (\triangle u)\cdot\,\na u \,dx \equiv \,
-\,\int_{\R^n_+}\, \pa_k (\pa^2_i\,u_j)\,\pa_k u_j\,dx=\,
\int_{\R^n_+}\,|\na^2\,u |^2\,dx -\,I
$$
where
$$
I\equiv \,\int_{\Ga} \,(\pa_i \pa_k \,u_j\,)\,(\pa_k \,u_j\,)\,\cn_i
\,d\Ga =\,-\,\int_{\Ga} \,(\pa_k \pa_n \,u_j\,)\,(\pa_k \,u_j\,)
\,d\Ga \,,
$$
since $\,\cn\,,$ the unit external normal to $\,\Ga\,,$ has
components $\,(0,\,...,\,0,\,-1)\,.$ If $j<n$ and $k=n$ the terms
$\pa_k \,u_j\,$ vanishes, due to the boundary conditions
\eqref{bboundary}. If $j<n\,,$ but $k<n\,,$ the terms $\,\pa_k \pa_n
\,u_j\,\,$ vanishes, since $\,\pa_n \,u_j=\,0\,$ on the boundary,
and $\,\pa_k $ is a tangential derivative. Hence we merely have to
consider the $\,j=\,n\,$ terms, namely $\,(\pa_n \pa_k \,u_n)\,\pa_k
\,u_n \,.$ If $k<n\,,$ it follows $\pa_k \,u_n =\,0\,.$ On the other
hand, if $\,k=\,n\,,$ by appealing to the divergence free condition,
one has
\begin{equation}%
\label{deln}%
\pa_n \pa_n \,u_n\,=\,-\,\sum_{j<\,n}\,\pa_n\,(\pa_j \,u_j\,)=\,0\,,
\end{equation}
since $\,\pa_j (\pa_n \,u_j)=\,0\,$ on $\,\Ga\,,$ for $\,j<\,n\,.$
We have shown that
$$
-\,\m\,\int_{\R^n_+}\, \na (\triangle u)\cdot\,\na u \,dx=\,\m\,
\int_{\R^n_+}\,|\na^2\,u |^2\,dx\,.
$$
the boundary integral related to the viscous term vanishes.\par%
Next we consider the pressure term. One has, by an integration by
parts,
$$
\int_{\R^n_+}\, (\,\na(\na \pi)\,)\cdot\,\na u \,dx \equiv \,
\int_{\R^n_+}\, \pa_k (\pa_j\,\pi)\,\pa_k u_j\,dx=\, -
\int_{\R^n_+}\,(\na \pi)\cdot\,\na\,(\na \cdot\,u\,) \,dx +\,A,
$$
where
$$
A \equiv \,\int_{\Ga} \,(\pa_k\,\pi)\,(\pa_k\,u_j) \,\cn_j
\,d\Ga=\,-\,\int_{\Ga} \,(\pa_k\,\pi)\,\pa_k\,u_n \,\,d\Ga =
\,-\,\int_{\Ga} \,(\pa_n\,\pi)\,(\pa_n\,u_n) \,d\Ga \,,
$$
since $\,\pa_k\,u_n =\,0\,$ on the boundary for $\,k<\,n\,.$
Furthermore, the volume integral on the right hand side vanishes,
due to the divergence free condition.\par%
Let's see that $\,A=\,0\,$ by showing that $\,\pa_n\,\pi=\,0\,$ on
$\,\Ga\,.\,$ By appealing to the n.th equation \eqref{nse-2} we show
that $\,\pa_n\,\pi =\,- \pa_t\,u_n - (u\cdot\nabla)\,u_n
+\,\m\,\triangle u_n\,.$ So, by appealing to boundary condition
$\,u_n=\,0\,,$ one easily shows that
$$
\pa_n\,\pi=\,\m\,\triangle \,u_n \quad \textrm {on} \,\Ga\,.
$$
Note that $\,(u\cdot\nabla)\,u_n =\,0\,$ on $\Ga\,.$ By taking into
account that the second order tangential derivatives of $\,u_n\,$
vanish on the boundary, we show that $\,\triangle \,u_n =\,0\,,$ by
appealing to \eqref{deln}. So $\,\pa_n\,\pi=\,0\,$ on $\,\Ga\,,$ as
desired. We have shown that
\begin{equation}
\label{presaa}%
\int_{\R^n_+}\, (\,\na(\na \pi)\,)\cdot\,\na u \,dx=\,0\,.
\end{equation}
Equation \eqref{estimum-2} is proved.\par%
Note that equation \eqref{presaa} holds under the non-slip boundary
condition, with a simpler proof. In fact, in this case, $\,A=\,0\,$
follows immediately from $\,\pa_n\,u_n =\,0\,$ on $ \,\Ga\,,$ which
is an immediate consequence of the divergence free property and the
non-slip boundary assumption.

\vspace{0.23cm}

The next, and main, step is to consider the non-linear term. We
start by showing that
\begin{equation}
\label{puns}%
\int_{\R^n_+} \,\na \,[(u\cdot\nabla)\,u\,]\cdot\na\,u \, dx
=\,\int_{\R^n_+} \,(\pa_k\,u_i)(\pa_i\,u_j)(\pa_k\,u_j)\,dx\,.
\end{equation}
This follows from the identity
\begin{equation}
\label{dvac}%
\,\na \,[(u\cdot\nabla)\,u\,]\cdot\na\,u
=\,(\pa_k\,u_i)(\pa_i\,u_j)(\pa_k\,u_j)+\, \frac12 \,u_i\,\pa_i
\,\big( \,\sum_{j,\,k} (\,\pa_k u_j)^2\,\big)
\end{equation}
since, by an integration by parts, we show that the integral of the
second term on the right hand side of the \eqref{dvac} vanishes, as
follows from the divergence free and the tangential to the boundary
properties (unfortunately, in the cylindrical coordinates case, the
counterpart of this main point is much more involved).

\vspace{0.2cm}

Next we prove the main estimate \eqref{bacho-2}. Following
\cite{baechoe}, we consider separately the three cases $\,i
\neq\,n\,;$ $\,i=\,n\,$ and $\,j \neq\,n\,;$ $\,i=\,j=\,n\,.$\par%
If $\,i \neq\,n\,,$ one has
\begin{equation}
\label{ineqn}%
\begin{aligned}
\int_{\R^n_+} \,(\pa_k\,u_i)(\pa_i\,u_j)(\pa_k\,u_j)\,dx=\\
-\,\int_{\R^n_+} \,u_i\,\pa_k\,\big(
(\pa_k\,u_j)(\pa_i\,u_j)\,\big)\,dx +\,\int_{\Ga} \,\,u_i\,
(\pa_k\,u_j)(\pa_i\,u_j)\,\nu_k \,dx
\end{aligned}
\end{equation}
The boundary integral is equal to
$$
-\,\int_{\Ga} \,\,u_i\, (\pa_n\,u_j)(\pa_i\,u_j) \,dx\,.
$$
If $\,j \neq\,n\,,$ one has $\,\pa_n\,u_j=\,0\,.$ If $\,j=\,n\,,$
one has $\,\pa_i\,u_j=\,0\,,$ since $\,\pa_i\,$ is a tangential
derivative and $\,u_n=0\,.$ Hence the boundary integral in equation
\eqref{ineqn} vanishes. On the other hand, since $\,i \neq\,n\,,$
the volume integral on the right hand side of equation \eqref{ineqn}
is bounded by the right hand side of inequality \eqref{bacho-2}.
After all, if $\,i \neq\,n\,,$ the left hand side of equation
\eqref{ineqn} is
bounded by the right hand side of inequality \eqref{bacho-2}.\par%
Next we assume that $\,i=\,n\,$ and $\,j \neq\,n\,.$ In this case,
by an integration by parts, one gets
\begin{equation}
\label{ienjn}%
\begin{aligned}
&\int_{\R^n_+} \,(\pa_k\,u_i)(\pa_i\,u_j)(\pa_k\,u_j)\,dx=\\
&-\,\int_{\R^n_+} \,(\triangle u)(\pa_i\,u_j)\,u_j \,dx\\
&-\,\int_{\R^n_+} \,(\pa_k\,u_n)(\pa_i\,\pa_k \,u_j) \,u_j \,dx\,,
\end{aligned}
\end{equation}
since the boundary integral, which appears after the above
integration by parts, vanishes. In fact, the terms
$\,(\pa_i\,u_j)\,$ vanish on the boundary, for  $\,i=\,n\,$ and $\,j
\neq\,n\,.$ From \eqref{ienjn} it follows that the left hand side of
this equation is bounded by  the right hand side of inequality
\eqref{bacho-2},
as desired.\par%
If  $\,i=\,j=\,n\,,$ we have to estimate the integral
$$
B \equiv \,\int_{\R^n_+} \,(\pa_k\,u_n)^2(\pa_n\,u_n)\,dx=\,
-\,\int_{\R^n_+} \,(\pa_k\,u_n)^2( \sum_{j\neq\,n}
\pa_j\,u_j)\,dx\,.
$$
By integration by parts one gets
$$
B=\,2\,\int_{\R^n_+} \,(\pa_k\,u_n)( \sum_{j\neq\,n}
\,\pa_j\,\pa_k\,u_n)\, u_j\,dx - \int_{\Ga} \,\,(\pa_k\,u_n)^2
\sum_{j\neq\,n}\,u_j\,\nu_j \,d\Ga\,.
$$
Since the above boundary integral vanishes, the absolute value of
$B$ is  bounded by  the right hand side of inequality
\eqref{bacho-2}. The proof of \eqref{bacho-2} is accomplished. From
now on the proof of theorem \ref{theocas} follows a very classical
way. For the readers' convenience we recall how to prove
\eqref{estimun}. From \eqref{estimum-2} and \eqref{bacho-2} it
follows that
\begin{equation}
\label{estibom}%
\frac12\,\frac{d}{dt}\,\int_{\R^n_+}\,|\na\,u|^2 \,dx +\,\m\,
\int_{\R^n_+}\,|\na^2\,u |^2\,dx \leq\, c(n)\,
\|\,|\ou|\,\na\,u\,\|_2 \, \|\,\na^2\,u \,\|_2\,.
\end{equation}
On the other hand, by H\H older's inequality,
$$
\|\,|\ou|\,\na\,u\,\|_2\leq\, \|\,\ou\,\|_p \,
\|\,\na\,u\,\|_{\frac{2\,p}{p-\,2}}\,.
$$
Furthermore, by interpolation and Sobolev's embedding theorem,
$$
\|\,\na\,u\,\|_{\frac{2\,p}{p-\,2}}\leq\,
\|\,\na\,u\,\|^{1-\frac{n}{p}}_2\,\|\,\na\,u\,\|^{\frac{n}{p}}_{2^*}\leq\,c\,
\|\,\na\,u\,\|^{1-\frac{n}{p}}_2\,\|\,\na^2\,u\,\|^{\frac{n}{p}}_2\,,
$$
since $\,(p-\,2)/(2\,p)=\,(1-\,n/p)/\,2 +\,(n/p)/2^*\,$. Here $\,2^*
=\,2\,n/(n-\,2)\,$ is a well known Sobolev's embedding exponent
(note that each single component of the tensor $\,\na\,u\,$
satisfies an homogeneous, Dirichlet or Neumann, boundary condition
on $\,\Ga $). Consequently,
$$
\|\,|\ou\,|\,\na\,u\,\|_2 \, \|\,\na^2\,u \,\|_2 \leq\,c\,
\|\,\ou\,\|_p
\,\|\,\na\,u\,\|^{1-\frac{n}{p}}_2\,\|\,\na^2\,u\,\|^{1+\,\frac{n}{p}}_2
\,.
$$

Hence, by Young's inequality,
\begin{equation}
\label{bombom}%
\|\,|\ou\,|\,\na\,u\,\|_2 \, \|\,\na^2\,u \,\|_2
\leq\,c\,\leq\,c\,\|\,\ou\,\|^q_p\,\|\,\na\,u\,\|^2_2\,+\,(\m
/2)\,\|\,\na^2\,u\,\|^2_2\,.
\end{equation}
From \eqref{estibom} and \eqref{bombom} we get, for
$\,t\in\,(0,\,T]\,,$
\begin{equation}
\label{estibom}%
\frac12\,\frac{d}{dt}\,\|\na\,u\|^2_2 +\,\frac{\m}{2}\,\|\na^2\,u
\|^2_2\,\leq\, c\,\|\,\ou\,\|^q_p\,\|\,\na\,u\,\|^2_2\,.
\end{equation}
This estimate immediately leads to \eqref{estimun} since, by the
Prodi-Serrin's assumption,
$$
\|\,\ou\,\|^q_p \in\,L^1(0,\,T)\,.
$$
\section{Non flat boundaries.}\label{tris}
It would be of basic interest to understand how crucial is in the
theorem \ref{theocas} the flat-boundary hypothesis. Does the result
hold in the neighborhood of non-flat boundary points? We propose to
start from the following particular case.
\begin{open problem}
\label{p1}%
\rm{ Let $\,(\rho,\,\theta,\,z)\,$ be the canonical cylindrical
coordinates in the three dimensional space, and consider the subset
defined by imposing to $\,\rho\,$ the constraint $\,
r<\,\rho<\,R\,,$ where $r$ and $R$ are positive constants. Assume
the slip boundary condition on the two lateral cylindrical surfaces,
and space periodicity with respect to the axial $\,z-$ direction. To
\emph{prove} or \emph{disprove} that Prodi-Serrin's condition on the
two "tangential" components $\,u_{\theta}\,$ and $\,u_z\,$ implies
smoothness.}
\end{open problem}
A careful, preliminary, study of the above problem showed that the
possible main obstacles still appear in a $\,2-D\,$ approach in
plane-polar coordinates. Furthermore, calculations in this direction
led us to be convinced that the reply to the open problem \ref{p1}
is positive.
\section{The limit case $\,p=\,n\,.$}\label{quatris}
In the following $\,\|\cdot\,\|_p\,$ denotes the $\,L^p\,$ norm,
moreover $\,\|\cdot\,\|=\,\|\cdot\,\|_2\,.$ The Prodi-Serrin's
condition for $\,(q,\,p)=\,(\infty,\,n)$, namely
\begin{equation}
\label{pslimo}%
u \in\, L^{\infty}(0,\,T;\,L^n(\Om)\,)\,,
\end{equation}
deserves a separate treatment. In reference \cite{sohr} uniqueness
of solutions was proved under assumption \eqref{pslimo}. In
references \cite{giga} and \cite{wahl} strong regularity was proved
by assuming time-continuity in $\,L^n(\Om)\,.$\par%
In reference \cite{bdv97} the system \eqref{nse-2} was considered in
a smooth bounded domain $\,\Om\,,$ under the \emph{non-slip}
boundary condition. It was shown that, under assumption
\eqref{pslimo}, the addition of a sufficiently small upper-bound for
possible left-discontinuities on the norm $\,\|u(t)\,\|_n\,$ implies
regularity (the same result was proved in the same year, in
reference \cite{sohr-2}). Note that, under assumption
\eqref{pslimo}, weak solutions are weakly continuous with values in
$\,L^n(\Om)\,.$ Hence the pointwise values $\, \|u(t)\|_n \,$ are
everywhere well defined. In reference \cite{bdv97} a very simple
explicit upper bound for the required left-discontinuities was also
shown. It was proved that there is a positive, independent, constant
$\,C\,$ such that if \eqref{pslimo} holds and, in addition,
\begin{equation}
\label{cces}%
\limsup_{t \rightarrow\,\ot -\,0} \|u(t)\|_n <\,\|u(\ot)\|_n
+\,4\,(C/4)^n
\end{equation}
for each $\,\ot \in \,(0,\,T]\,,$ then solutions are smooth. Let's
exhibit a particularly explicit value $\,C\,.$ Consider the main
case $\,n=\,3\,$. Denote by $\,A\,$ the well known operator
$\,A=\,-\,P\,\De\,$, where $\,P\,$ is the classical orthogonal
projection of $\,L^2\,$ onto $\,H\,,$ the closure in $\,L^2\,$ of
smooth, compact supported functions in $\,\Om\,.$ Let $\,c_0\,$ be
such that
$$
\|\,\na\,v\,\|_6 \leq\,c_0\,\|A\,v\|\,, \quad \forall \,v
\in\,D(A)\,.
$$
Note that, roughly speaking, the constant $c_0$ is like the constant
$\,c_1\,$ in the Sobolev's embedding $\,\|\,v\,\|_6 \leq\,c_1\,
\|\na\,\,v\|\,.$ We have shown that the constant
\begin{equation}
\label{mahh}%
C=\,\frac{\m}{2\,c_0}\,,
\end{equation}
in equation \eqref{cces}, implies smoothness.\par%
More recently, in the famous article \cite{seregin}, the authors
succeed in proving that, in the whole space case, the assumption
\eqref{pslimo} alone guarantees smoothness of solutions. This was,
for many years, one of the most challenging, and difficult, open
problems in the mathematical theory of Navier-Stokes equations.
Later on, extensions to the boundary has been obtained. See
\cite{seregin-2} for the half-space case, and \cite{shilkin} for
curved smooth boundaries. Proofs of the above results are quite
involved. So, despite these strong improvements, a look at reference
\cite{bdv97} could be of some interest, since the straightforward
proof is particularly elementary.

\vspace{0.2cm}

Concerning the subject of the present notes, namely statements
involving only $\,n-1\,$ components of the velocity, we recall that
the results proved in \cite{bdv97} were extended in reference
\cite{bdv2000} to the case in which the two conditions,
\eqref{pslimo} and \eqref{cces}, are merely required to $\,\ou\,.$
It would be not difficult to extend to the half-space case, under
slip boundary conditions, the result stated in \cite{bdv2000}, by
putting together the proof developed in this last reference and that
shown in the present notes.

\vspace{0.2cm}

It would be interesting to study the following problem, even in the
whole space case.
\begin{open problem}
\label{p2}%
\rm{ To extend the main result proved in reference \cite{seregin} to
the "two components case".}
\end{open problem}

\end{document}